\documentclass[12pt]{article}
\usepackage{amsfonts}
\usepackage[dvips]{graphicx}
\usepackage{amssymb,color}
\usepackage[latin1]{inputenc}
\usepackage{epic}
\usepackage[latin1]{inputenc}
\usepackage{enumerate}
\usepackage{color}
\usepackage{lettrine}
\usepackage{calc,pifont}
\usepackage[noindentafter,calcwidth]{titlesec}
\usepackage{amsmath}
\usepackage{amssymb}
\usepackage{amsthm}
\usepackage{subfigure}

\def\BBox{\kern  -0.2cm\hbox{\vrule width 0.2cm height 0.2cm}}

\newtheorem{theorem}{Theorem}[section]

\newtheorem{corollary}{Corollary}[section]

\newtheorem{conjecture}[theorem]{Conjecture}

\begin{document}

\title{A note on Erdös-Faber-Lovász Conjecture and edge coloring of complete graphs.
}
\author{
G. Araujo-Pardo\thanks{garaujo@math.unam.mx}, A. Vázquez-Ávila\thanks{pare\_23@hotmail.com}\\{\small  Instituto de Matem\'aticas}\\ {\small Universidad Nacional
Auton\'oma de M\'exico}\\
}

\date{}
\maketitle
\begin{abstract}
A linear hypergraph is intersecting if any two different edges have
exactly one common vertex and an $n$-quasicluster is an intersecting
linear hypergraph with $n$ edges each one containing at most $n$
vertices and every vertex is contained in at least two edges. The
Erdös-Faber-Lovász Conjecture states that the chromatic number of
any $n$-quasicluster is at most $n$. In the present note we prove
the correct\-ness of the conjecture for a new infinite class of
$n$-quasiclusters using a specific edge coloring of the complete
graph.
\end{abstract}




\section*{Introduction}
A \emph{hypergraph} $\mathcal{H}=(\mathcal{V},\mathcal{E})$ consists
of a finite no empty set $\mathcal{V}$, the vertices of
$\mathcal{H}$, and a finite collection  $\mathcal{E}$ of subsets of
$\mathcal{V}$, the edges of $\mathcal{H}$. It is assumed that each
vertex belongs to at least one edge. A hypergraph
$\mathcal{H}=(\mathcal{V},\mathcal{E})$ is \emph{linear} if $|E\cap
F|\leq1$, for all $E,F\in \mathcal{E}$, and $\mathcal{H}$ is
\emph{intersecting} if $|E\cap F|=1$, for all $E,F\in \mathcal{E}$.
In the remainder of this note each edge contains at least two
vertices.

Let $\{1,2,\ldots,k\}$ be a set of $k$ colors. A
$k$\emph{-vertex-coloring} of $\mathcal{H}$ is a surjective map
$\varphi:\mathcal{V}\longrightarrow\{1,2,\ldots,k\}$ such that if
$u,v\in\mathcal{V}$ are adjacent, then $\varphi(u)\neq\varphi(v)$.
In other words, no two vertices with the same color belong in the
same edge. \emph{The chromatic number} of $\mathcal{H}$, denoted by
$\chi(\mathcal{H})$, is the minimum $k\in\mathbb{N}$ for which there
is a $k$-vertex-coloring of $\mathcal{H}$.

Thus, in the hypergraph setting the original Erdös-Faber-Lovász
Conjecture states:

\begin{conjecture}\label{original}
If $\mathcal{H}$ is a linear hypergraph consisting of $n$ edges,
each one containing $n$ vertices, then $\chi(\mathcal{H})=n$.
\end{conjecture}

An $n$-\emph{quasicluster} is an intersecting linear hypergraph
consisting of $n$ edges, each one with at most $n$ vertices and each
vertex is contained in at least two edges.

\begin{conjecture}\label{cuasicumulo}
If $\mathcal{H}$ is an $n$-quasicluster, then $\chi(\mathcal{H})\leq
n$.
\end{conjecture}

It is not difficult to prove that the conjecture \ref{original} and
the conjecture \ref{cuasicumulo} are equivalent. For this, consider
an $n$-quasicluster and add new vertices in each edge such that the
edges of this new intersecting linear hypergraph $\mathcal{H}$ has
size $n$. By assumption Conjecture \ref{original} is true for
intersecting linear hypergraphs, the hypergraph $\mathcal{H}$ has an
$n$-vertex-coloring, which in turn induces an $n$-vertex-coloring to
$n$-quasicluster. On the other hand, consider a hypergraph
$\mathcal{H}$ which satisfies the conjecture \ref{original}. By the
theorem $3$ of \cite{Romero2} there exists an intersecting linear
hypergraph $\mathcal{\widetilde{H}}$ consisting of $n$ edges, each
of size $n$, such that
$\chi(\mathcal{H})=\chi(\mathcal{\widetilde{H}})$. Now, deleting
from $\mathcal{\widetilde{H}}$ vertices of degree one we get an
$n$-quasicluster $\mathcal{H}'$. Assuming that Conjecture
\ref{cuasicumulo} is true we have $\chi(\mathcal{H}')\leq n$ and
this coloring can be easily extended to $n$-vertex-coloring to
$\mathcal{\widetilde{H}}$ (using non-used colors to each vertex of
degree one), thus $\chi(\mathcal{H})=n$.

In the present note we say that $\mathcal{H}$ is an \emph{instance}
of the conjecture or theorem, if $\mathcal{H}$ satisfies the
hypothesis of its statements.


There exist works related with some equivalences of Conjecture
\ref{original} and also many advances, but it is clear that its
proof is, in this moment, far from being attained. There are some
results about upper bonds on the number of colors required.
Specifically, Mitchem \cite{Mitchem 2}, and independently Chang and
Lawler \cite{Chang} had shown that if $\mathcal{H}$ is an instance
of the conjecture \ref{original}, then the chromatic number of
$\mathcal{H}$ is at most $\lceil\frac{3n}{2}-2\rceil$. Kahn
\cite{Kahn} had proved, as an asymptotic result, that if
$\mathcal{H}$ is an instance of the conjecture \ref{original}, then
the chromatic number of $\mathcal{H}$ is at most $n+o(n)$. There are
some works about specific classes of hypergraphs that satisfy the
conjecture \ref{original}; see for example \cite{C. Berge},
\cite{Colobourn}, \cite{Jackson}, \cite{Mitchem}, \cite{Romero2} and
\cite{arroyo}. Also there are interesting equivalences of this
conjecture; see for example  \cite{haddad}, \cite{Hindman},
\cite{margraf} and \cite{Romero1}. Recently, Faber \cite{Faber}
proved that for regular and uniform linear hypergraphs of fixed
degree there can only be a finite number of counterexamples for
conjecture \ref{original}.

In this work we expose a new method to approach it, using a specific
edge coloring of the complete graph, giving a new infinite class of
$n$-quasiclusters that satisfy the conjecture \ref{cuasicumulo}.
\section{The result}

Let $G$ be a simple graph. A \emph{decomposition} of $G$ is a
collection $\mathcal{D}=\{G_{1},\ldots,G_{k}\}$ of subgraphs of $G$
such that every edge of $G$ belongs to exactly one subgraph in
$\mathcal{D}$, denoted by $(G,\mathcal{D})$ a decomposition of $G$.

Let $\{1,\ldots,k\}$ be a set of $k$ colors. A
$k$-$\mathcal{D}$-\emph{coloring} of $(G,\mathcal{D})$ is a
surjective map $\varphi':\mathcal{D}\longrightarrow\{1,\ldots,k\}$
such that for every $G,H\in\mathcal{D}$ if $V(G)\cap
V(H)\neq\emptyset$ then $\varphi'(G)\neq \varphi'(H)$. Here
$\varphi'$ means that every edge of the subgraph $G$ is colored with
the color $\varphi'(G)$. The \emph{chromatic index} of a
decomposition $(G,\mathcal{D})$, denoted by
$\chi'((G,\mathcal{D}))$, is the minimum $k\in\mathbb{N}$ for which
there is a $k$-$\mathcal{D}$-coloring of $(G,\mathcal{D})$. Now the
$(K_{n},\mathcal{D})$ denotes a decomposition where the elements of
$\mathcal{D}$ are complete subgraphs of $K_{n}$.

It is not difficult to see that there exists a bijection between the
decompositions $(K_{n},\mathcal{D})$ and the $n$-quasiclusters
\footnote{The elements of the decomposition $(K_n,D)$ can be thought
as the cliques of the intersection graph of the corresponding
$n$-quasicluster (see \cite{Hindman}, \cite{margraf},
\cite{Mitchem}, \cite{wiki}).}. Because, if we have a decomposition
$(K_{n},\mathcal{D})$, each vertex of $K_n$ is associated with an
edge of the $n$-quasiclusters and each element $G\in\mathcal{D}$
with a vertex of the $n$-quasiclusters; the intersection vertex of
the edges associated with the vertices of $G$. Then a
$k$-vertex-coloring of an $n$-quasicluster $\mathcal{H}$ is a
$k$-$\mathcal{D}$-\emph{coloring} of $(K_{n},\mathcal{D})$ and
vice-versa. Therefore, the  Conjecture \ref{cuasicumulo} is
equivalent to:

\begin{conjecture}
If $(K_{n},\mathcal{D})$ is a decomposition, then
$\chi'((K_{n},\mathcal{D}))\leq n$.
\end{conjecture}
As example of our interpretation consider an $n$-quasicluster where
each vertex is a member of exactly two edges, that is, each vertex
has degree two, then the elements of corresponding decomposition
$(K_{n},\mathcal{D})$ are subgraphs of order two, namely
$\mathcal{D}=E(K_n)$. Then a
$\chi'((K_{n},\mathcal{D}))$-edge-coloring of $(K_{n},\mathcal{D})$
is equivalent to $\chi'(K_{n})$-edge-coloring of $K_{n}$ and by
Vizing's Theo\-rem we have that $\chi'((K_{n},\mathcal{D}))\leq n$.

In this note we will work with a specific $n$-edge-coloring of the
complete graph $K_{n}$ given by the following: suppose that
$K_{n}=(\mathbb{Z}_{n},E)$ and that $\{c_{1},\ldots,c_{n}\}$ is a
set of $n$ different colors. If $ab\in E$ is an edge then the
associated color for this edge is $c_{a+b}$ (with $a+b\in
\mathbb{Z}_n$). This assignment is an $n$-edge-coloring of $K_n$
(for every $n\geq2$). Let $G_{0},\ldots,G_{n-1}$ be the $n$
chromatic classes of $K_n$ with respect to this coloring. If we
think that $G_{i}=(\mathbb{Z}_{n},E_{i})$, where $E_{i}=\{ab\in
E:\quad\!\!\!a+b\equiv i\quad\!\!\!\mbox{mod $n$}\}$ for
$i=0,\ldots,n-1$, then these subgraphs satisfy:

\begin{figure}[]
\begin{center}
\includegraphics[height =6cm]{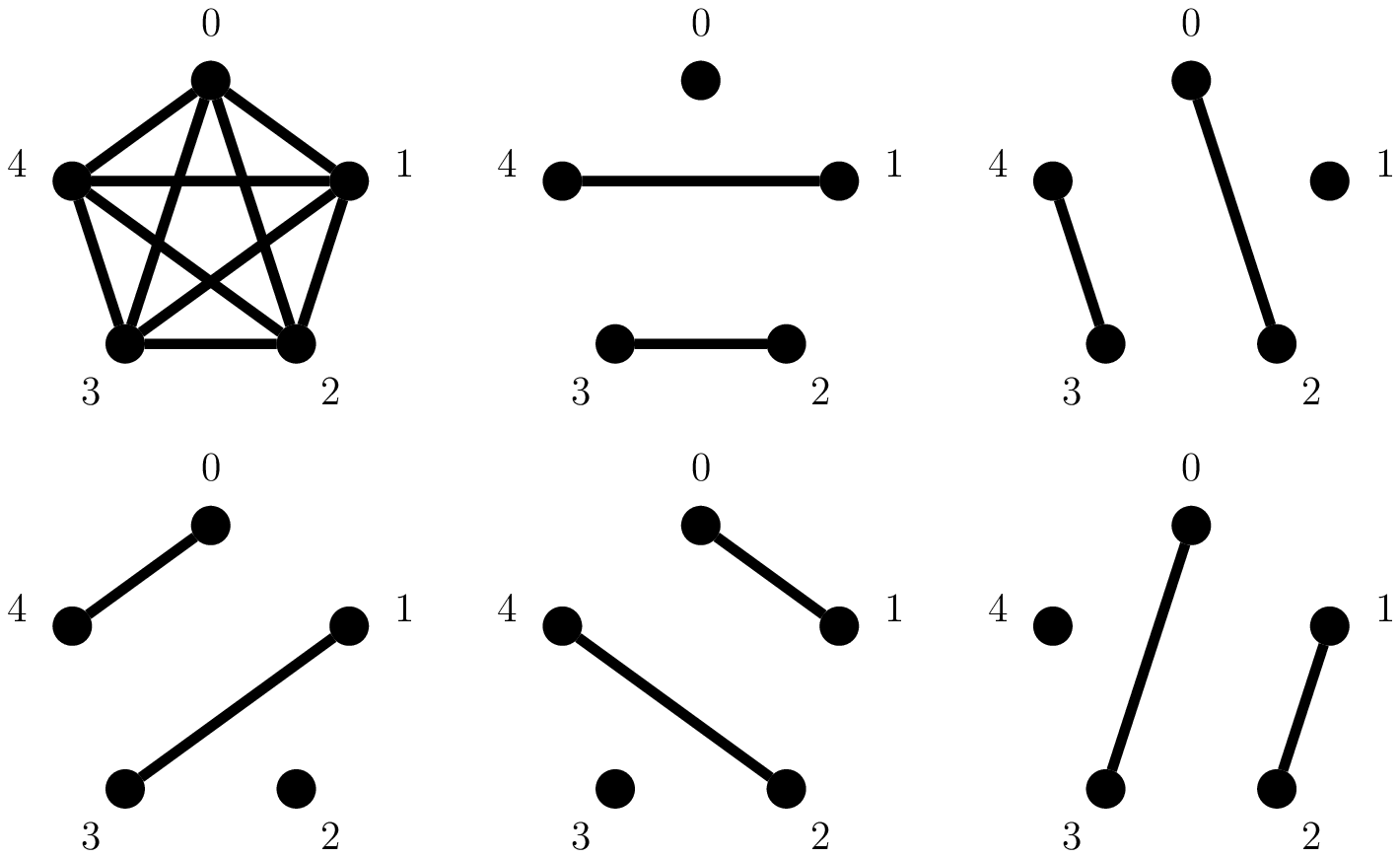}
\caption{Subgraphs corresponding to the $5$-edge-coloring of
$K_{5}$.}\label{des5}
\end{center}
\begin{center}
\includegraphics[height =5cm]{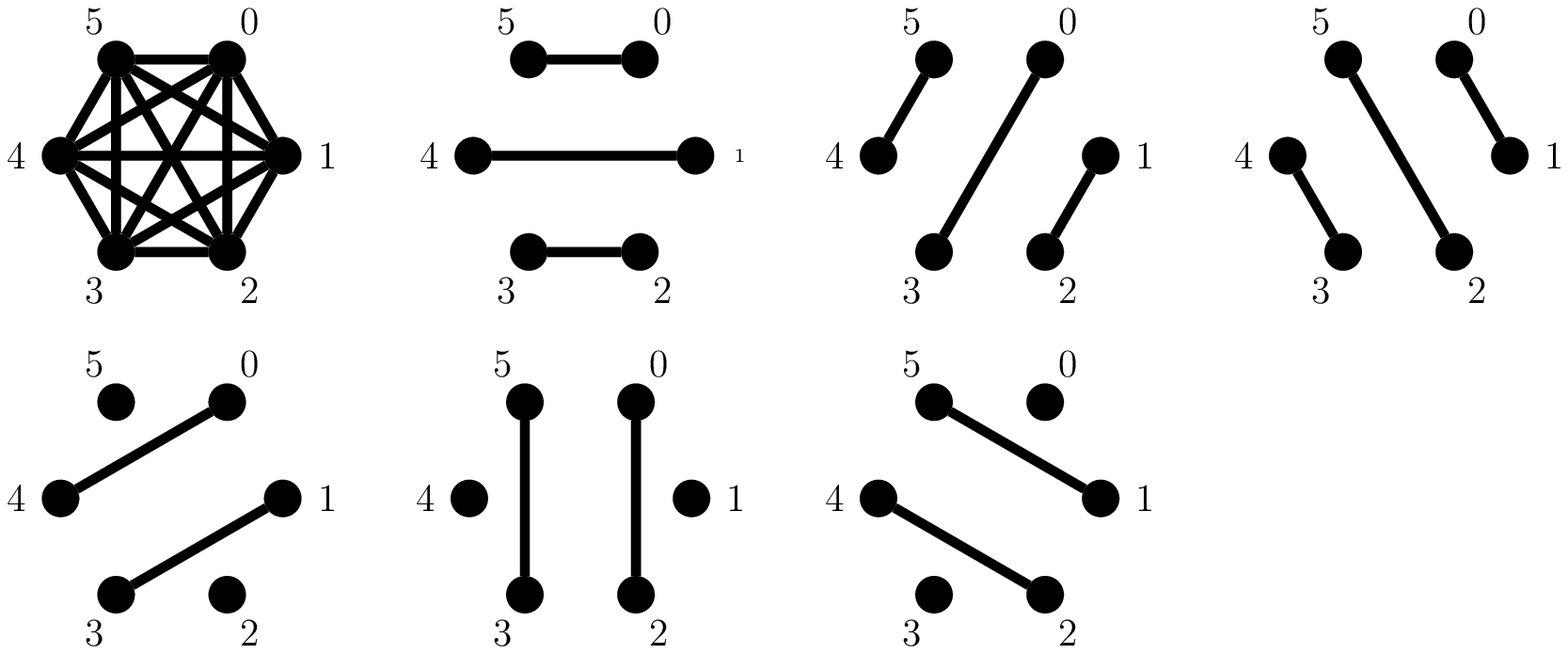}
\caption{Subgraphs corresponding to the $6$-edge-coloring of
$K_{6}$.}\label{des6}
\end{center}
\end{figure}

\begin{enumerate}
 \item For $i=0,\ldots,n-1$, the degree of the vertices of $G_{i}$ is at most one.
 \item If $n$ is odd then the subgraph $G_{i}$ has an isolated vertex,
 say $u_{i}$ and $G_{i}-u_{i}$ is a perfect matching. If
 $n$ is even and $i$ is even then the subgraph $G_{i}$ has two
 isolated vertices, say $u_{i}$ and $v_{i}$ and $G_{i}-u_{i}-v_{i}$ is a perfect matching. When $i$ is odd the subgraph $G_{i}$
 is a perfect matching.
 \end{enumerate}
For example, Figures \ref{des5} and \ref{des6} show the subgraphs
corresponding to $K_5$ and $K_6$, res\-pectively.

Our results are related with a previous result given by Romero and
Sánchez-Arroyo in \cite{Romero2}. We define some similar concepts
like them, but their approach to solve this is totally algorithmic.
On the other hand, our method, as we previously exposed, is related
with edge colorings in the complete graph.

In \cite{Romero2} the following concepts were defined: A nonempty
set $W$ of nonnegative integers is \emph{compact}, if either
$|W|=1$, or there is an order $(a_{1},\ldots,a_{|W|})$ on $W$, such
that $a_{i+1}=a_{i}+1$, for $i=1,\ldots,|W|-1$. Let
$\mathcal{H}=(\mathcal{V},\mathcal{E})$ be a linear hypergraph with
$n$ edges. $\mathcal{H}$ is \emph{edge conformable} if there is a
bijection $\psi:\mathcal{E}\longrightarrow\{0,\ldots,n-1\}$, called
\emph{conformal labeling}, such that for each vertex
$v\in\mathcal{V}$, the set $F(v)=\{\psi(E):v\in E\in \mathcal{E}\}$
can be partitioned into two compact sets. The main result of
\cite{Romero2} is that any intersecting li\-near hypergraph
consisting of $n$ edges, each of size $n$, and edge conformable has
an $n$-vertex-coloring.


Now we introduce the following definitions: Let
$W=\{w_{1},\ldots,w_{r}\}$ be a subset of $\mathbb{Z}_{n}$;  $W$ is
$k$-\emph{arithmetic} if $w_{i+1}-w_{i}\equiv k$ mod $n$, for
$i=1,\ldots,r-1$ and some $k\in
\{1,\ldots,\lfloor\frac{n}{2}\rfloor\}$. A decomposition
$(K_{n},\mathcal{D})$ is called \emph{arithmetic decomposition}, if
there exists a bijection
$\varphi:V(K_{n})\longrightarrow\mathbb{Z}_{n}$, called
\emph{arithmetic labeling}, such that for every $G\in\mathcal{D}$
either $V(G)$ is $k$-arithmetic or $V(G)$ can be partitioned into
two $k$-arithmetics sets of same cardinality, for some
$k\in\{1,\ldots,\lfloor\frac{n}{2}\rfloor\}$. Let
$(K_{n},\mathcal{D})$ be an arithmetic decomposition and $G$ be an
element of $\mathcal{D}$, where $V(G)=\{v_{1},\ldots,v_{l}\}$ has
odd cardinality, then $v_{\frac{l+1}{2}}$ is called the
\emph{central vertex} of $G$. Also, we say that a decomposition
$(K_n,\mathcal{D})$ has different central vertices if any pair of
central vertices (corresponding to elements of $\mathcal{D}$ of odd
order) are different.


\begin{theorem}\label{teoaritmetico1}
Let $(K_{n},\mathcal{D})$ be an arithmetic decomposition with
different central vertices, then $\chi'((K_n,\mathcal{D}))\leq n$.
\end{theorem}

\textit{\textbf{Proof}}

Let $(K_{n},\mathcal{D})$ be an arithmetic decomposition, $G$ be an
element of $\mathcal{D}$ and $\{c_1\ldots,c_n\}$ be a set of $n$
different colors.

\begin{description}

\item[Case (i)]  If $V(G)=\{v_{1},\ldots,v_{l}\}$ is $k$-arithmetic, for some $k\in\{1,\ldots,\lfloor\frac{n}{2}\rfloor\}$ then by hypothesis $v_{i+1}-v_{i}\equiv k$ mod $n$, for $i=1,...,l-1$. The edges of $G$ will be colored as follows:

\begin{description}
  \item[(i.a)] If $V(G)$ has even cardinality, then $v_{1}+v_{l}\equiv v_{2}+v_{l-1}\equiv\ldots\equiv v_{\frac{l}{2}}+v_{\frac{l}{2}+1}\equiv j$, for some $j\in\mathbb{Z}_{n}$. Assign the color $c_j$ to the edges $M=\{v_{1}v_{l},v_{2}v_{l-1},\ldots,v_{\frac{l}{2}}v_{\frac{l}{2}+1}\}$. As $M$ is a perfect matching then there are no incident edges to $M$ of color $c_j$ different from $M$, so that we can assign the color $c_{j}$  to all edges of $G$.
\

  \item[(i.b)] If $V(G)$ has odd cardinality, let $v_G=v_{\frac{l+1}{2}}\in V(G)$ be the central vertex of $G$. Then $v_{1}+v_{l}\equiv v_{2}+v_{l-1}\equiv\ldots\equiv v_{\frac{l-1}{2}}+v_{\frac{l+1}{2}+1}\equiv j$, for some $j\in\mathbb{Z}_{n}$, so that we can assign the color $c_j$  to all edges of $G-v_{G}$ as Case (i.a).
  \end{description}

      Note that the color $c_{j}$ is not incident to $v_G$: otherwise there exists $u_{G}\in V(K_{n})$ such that $v_{G}+u_{G}\equiv j$ mod $n$. As $V(G)$ is $k$-arithmetic and $v_{G}$ is the central vertex of $G$ then $(v_{G}-rk)+(v_{G}+rk)\equiv j$ mod $n$, for $r=1,\ldots\frac{|V(G)|-1}{2}$, that is $2v_G\equiv j$ mod $n$. Since $2v_G\equiv j$ mod $n$ and by hypothesis $v_{G}+u_{G}\equiv j$ mod $n$ then $2v_G\equiv v_{G}+u_{G}$ mod $n$; this implies that $u_G\equiv v_G$ mod $n$, which is a contradiction. Therefore there are no incident edges to $v_{G}$ of color $c_j$, and so we can assign the color $c_j$  to all edges of $G$.

\item[Case (ii)] Now, if $V(G)$ can be partitioned into two $k$-arithmetic sets of same cardinality, for some $k\in\{1,\ldots,\lfloor\frac{n}{2}\rfloor\}$ then the subgraph $G$ will be colored as follows: suppose that $V(G)=\{v_{1},\ldots,v_{l}\}\cup\{u_{1},\ldots,u_{l}\}$. By hypothesis $\{v_{1},\ldots,v_{l}\}$ and $\{u_{1},\ldots,u_{l}\}$ are $k$-arithmetics, for some $k\in\{1,\ldots,\lfloor\frac{n}{2}\rfloor\}$. Then $v_{i+1}-v_{i}\equiv u_{i+1}-u_{i}\equiv k$ mod $n$, for $i=1\ldots,n-1$; so that $v_{1}+u_{l}\equiv v_{2}+u_{l-1}\equiv\ldots\equiv v_{l}+u_{1}\equiv j$, for some $j\in\mathbb{Z}_{n}$. Assign the color $c_j$ to the edges $M=\{v_{1}u_{l},v_{2}u_{l-1},\ldots,v_{l}u_{1}\}$. As $M$ is a perfect matching, there are no incident edges to $M$ of color $c_j$ different from $M$, so that we can assign the color $c_{j}$  to all edges of $G$.
\end{description}

Remains prove that for all $G,H\in\mathcal{D}$, with $V(G)\cap
V(H)=\{v\}$ and different central vertex (in case of having it) have
different colors. Let $G,H\in\mathcal{D}$, then
\begin{description}
  \item[Case (i)] If $V(G)$ and $V(H)$ have even cardinality then the corresponding perfect matching of $G$ and $H$ does not share edges (by linearity), therefore they have different colors.
  \item[Case (ii)] Suppose that $V(G)$ and $V(H)$ have odd cardinality and the edges of $G$ have the same color that the edges of $H$. Let $v_G$ and $v_H$ be the central vertices of $G$ and $H$ respectively.
      \begin{description}
        \item[(ii.a)] If $\{v\}=V(G)\cap V(H)\not\subset\{v_G,v_H\}$ then there exists $u_G\in V(G)$ and $u_H\in V(H)$ such that $u_G+v\equiv u_H+v\equiv j$ mod $n$, for some $j\in\mathbb{Z}_n$, implying that $u_G\equiv u_H$ mod $n$, which is a contradiction.
        \item[(ii.b)] If $v_H=V(G)\cap V(H)$ then there exists $u_G\in V(G)$ such that $u_G+v_H\equiv(v_{H}-rk)+(v_{H}+rk)$ mod $n$, for $r=1,\ldots\frac{|V(G)|-1}{2}$. Since $2v_H=(v_{H}-rk)+(v_{H}+rk)$, for $r=1,\ldots\frac{|V(G)|-1}{2}$ then $u_G\equiv v_H$ mod $n$, which is a contradiction.
      \end{description}
Therefore $E(G)$ and $E(H)$ have different colors.
  \item[Case (iii)] Finally, if $V(G)$ and $V(H)$ have different
  cardinality it is not difficult to see that the edges of $G$
  and $H$ have different colors because if $u_G\in V(G)$ is the
  central vertex of $G$ then the perfect matching of $G-u_G$
  and the perfect matching of $H$ does not share edges (by linearity),
  and so they have different colors.\qed
\end{description}

\begin{figure}[ht!]
\begin{center}
\includegraphics[height =5cm]{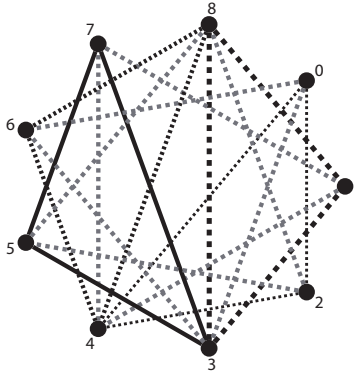}
\caption{Elements of $\mathcal{D}$ with order larger than
two.}\label{triangulos}
\end{center}
\end{figure}

To continue, in Figure \ref{triangulos}, we exhibit the theorem
giving an example. Let $V(G_{0})=\{0,3,6\}$, $V(G_{1})=\{1,4,7\}$,
$V(G_{2})=\{5,8,2\}$, $V(H_{0})=\{0,2,4\}$ $V(H_{1})=\{4,6,8\}$,
$V(H_{2})=\{8,1,3\}$ and $V(H_{3})=\{3,5,7\}$ be the vertices of the
complete graphs of $\mathcal{D}$ of cardinality larger than two and
the rest of the elements of $\mathcal{D}$ are edges.

Note that:
\begin{enumerate}
\item $V(G_{i})$ is $3$-arithmetic and $V(H_{j})$ is $2$-arithmetic,
for $i=0,1,2$ and $j=0,1,2,3$.
\item The central vertices of $V(G_{0})$, $V(G_{1})$, $V(G_{2})$, $V(H_{0})$,
$V(H_{1})$, $V(H_{2})$, $V(H_{3})$ are $3,4,8,2,6,1$ and $5$
respectively.
\end{enumerate}

Hence, this decomposition is a $(K_9,\mathcal{D})$ arithmetic
decomposition. By the Theorem \ref{teoaritmetico1} this
decomposition satisfies the conjecture \ref{cuasicumulo}.

To finish this note, it is important to establish which is the
correspondence, regarding the previous definitions, given in
arithmetic decompositions to  $n$- quasiclusters. We state the
Theorem \ref{teoaritmetico1} in these terms. To do this we give the
definitions in terms of hypergraphs (or $n$-quasiclusters).

Let $\mathcal{H}=(\mathcal{V},\mathcal{E})$ be an $n$-quasicluster.
$\mathcal{H}$ is \emph{edge arithmetic} if there is a bijection
$\varphi:\mathcal{E}\rightarrow\mathbb{Z}_n$, that we call
arithmetic labeling, such that for each vertex $u\in\mathcal{V}$,
the set $F(u)=\{\varphi(E):u\in E\in\mathcal{E}\}$ is $k$-arithmetic
or can be partitioned in two $k$-arithmetic sets of same
cardinality. Let $\mathcal{H}$ be an $n$-quasicluster edge
arithmetic, $u$ be a vertex of $\mathcal{H}$ of odd degree and
$F(u)=\{E_1,\ldots,E_l\}$, then $E_{\frac{l+1}{2}}$ is called
\emph{central edge} of $u$. We say that an $n$-quasicluster edge
arithmetic has different central edges if any pair of central edges
(corresponding to vertices of odd degree) is different. Then, the
main result (Theorem \ref{teoaritmetico1}) in terms of hypergraphs
states:

\begin{theorem}
Let $\mathcal{H}$ be an $n$-quasicluster. If $\mathcal{H}$ is edge
arithmetic and has different central edges, then
$\chi(\mathcal{H})\leq n$.
\end{theorem}

Finally we can note that any edge arithmetic $n$-quasicluster
$\mathcal{H}$ which has at most one vertex of odd degree in each
edge immediately has different central edges and then we have the
following:

\begin{corollary}
Let $\mathcal{H}$ be an $n$-quasicluster edge arithmetic with all
the edges with at most one vertex of odd degree, then
$\chi(\mathcal{H})\leq n$.
\end{corollary}
\subsection*{Acknowledgment}
The authors want to thank the anonymous referees for their kind help
and valuable suggestions which led to an improvement of this paper.

Also, the authors want to thank to CINNMA (Centro de Innovación
Matemática) and their people by academic support and friendship.

Research supported by PAPPIT-México under project IN101912



\begin{thebibliography}{\emph{Referencias}}\addcontentsline{toc}{section}{\textbf{\emph{Referencias}}}

\bibitem{Berge}C. Berge, \emph{Graphs}, Elsevier Science Publishers B.V.
North-Holland, $1991$.
\bibitem{C. Berge}C. Berge, A.J.W. Hilton, \emph{On two conjectures on edge colouring hypergraphs}, Congressus Numer. $\textbf{70}$
($1990$) 99-104.
\bibitem{Jackson} J. A. Bondy, U. S. R. Murty, \emph{Graph Theory with Applications}, The MacMillan Press Ltd., $1976$.

\bibitem{Colobourn} C.J. Colbourn, M.J. Colbourn, \emph{The chromatic
index of cyclic Steiner $2$-designs}, Internat. J. Math. \& Math
Sci. $\textbf{5}$/$\textbf{4}$ ($1982$) 823--825.

\bibitem{Mat} M.S. Colbourn, R.A Mathon, \emph{On cyclic steiner
$2$-design}, Annals of Discrete Mathematics $\textbf{7}$ ($1980$)
215--253.

\bibitem{Chang} W.I. Chang, E.L. Lawler, \emph{Edge coloring of hypergraphs and a conjecture of
Erdös, Faber, Lovász}, Combinatorica $\textbf{8}$ ($1988$) 293--295.

\bibitem{Chung} F. Chung, R. Graham, \emph{Erdös on Graphs: His Legacy of Unsolved Problems} (A.K. Peters, Wellesley, MA, 1988),
97--99.

\bibitem{favoritos} P. Erdös. \emph{On the combinatorial problems which I
would most like to see solved}, Combinatorica $\textbf{1}$ ($1981$)
313--318.
\bibitem{erdos} P. Erdös, \emph{Problems and results in graph theory and
combinatorial analysis}, Discrete Mathematics $\textbf{72}$ ($1988$)
81--92.

\bibitem{Faber} V. Faber, \emph{The Erdös-Faber-Lovász conjecture- the uniform regular case},
Journal of Combinatorics $\textbf{1}$ ($2010$) 113--120.

\bibitem{haddad} L. Haddad, C. Tardif, \emph{A clone-theoric formulation of the Erdös-Faber-Lovász Conjecture},
Discussiones Mathematicae Graph Theory $\textbf{24}$/$\textbf{3}$
(2004) 545--549.

\bibitem{Hindman}N. Hindman, \emph{On a conjecture of Erdös, Faber and
Lovász about $n$-colourings}, Canadian J. Math. $\textbf{33}$
($1981$) 545--549.

\bibitem{Jackson} B. Jackson, G. Sethuraman, C.Whitehead, \emph{A note on the Erdös-Faber-Lovász},
Discrete Mathematicas $\textbf{307}$ ($2007$) 911--915.

\bibitem{Kahn} J. Kahn, \emph{Coloring nearly-disjoint hypergraphs with $n+o(n)$ colors},
Journal of combinatorial theory, Series A $\textbf{59}$ ($1992$)
31--39.

\bibitem{margraf} H. Klein, M. Margraf,\emph{ A remark on the conjecture of
Erdös, Faber and Lovász}, Journal of Geometry $\textbf{88}$ (2008)
116--119.

\bibitem{Mitchem 2} J. Mitchem,\emph{On n-coloring certain finite set systems}, Ars Combinatoria $\textbf{5}$ (1978)
207--212.

\bibitem{Mitchem} J. Mitchem, R.L. Schmidt,\emph{ On the Erdös-Faber-Lovász Conjecture}, Ars Combinatoria $\textbf{97}$ (2010)
497--505.

\bibitem{Romero2}D. Romero and A. Sánchez-Arroyo, \emph{Adding evidence to
the Erdös-Faber-Lovász conjecture}, Ars Combinatoria $\textbf{85}$
($2007$) 71--84.
\bibitem{Romero1} D. Romero, A. Sánchez-Arroyo, \emph{Advances on the
Erdös-Faber-Lovász conjecture}, in G. Grimmet; C. McDiarmid,
Combinatorics, Complexity, and Chance: A Tribute to Dominic Welsh,
Oxford Lecture Series in Mathematics and its Applications, Oxford
University Press $\textbf{34}$ $2007$ 285--298.
\bibitem{arroyo} A. Sánchez-Arroyo, \emph{The Erdös-Faber-Lovász
conjecture for dense hypergraphs}, Discrete Mathematics
$\textbf{308}$ (2008) 991--992.
\bibitem{Seymor} P.D. Seymour, \emph{Packing nearly-disjoint sets},
Combinatorica \textbf{2} ($1982$) 91--97.

\

{\Large\textbf{Cyber Reference}}

\

\bibitem{wiki} \verb"http://en.wikipedia.org/wiki/Erd%C5%91s%E2%80%93Faber"\\\verb"%E2%80%93Lov%C3%A1sz_conjecture#CITEREFHindman1981"


\end{thebibliography}
\end{document}